# Selberg Integrals, Multiple Zeta Values and Feynman Diagrams

David H. Wohl

December 11, 2001

=== *Work in Progress* ===

BS"D

## ABSTRACT

We prove that there is an isomorphism between the Hopf Algebra of Feynman diagrams and the Hopf algebra corresponding to the Homogenous Multiple Zeta Value ring H in C<<X,Y>> . In other words, Feynman diagrams evaluate to Multiple Zeta Values in all cases. This proves a recent conjecture of Connes-Kreimer, and others including Broadhurst and Kontsevich.

The key step of our theorem is to present the Selberg integral as discussed in Terasoma [22] as a Functional from the Rooted Trees Operad to the Hopf algebra of Multiple Zeta Values. This is a new construction which provides illumination to the relations between zeta values, associators, Feynman diagrams and moduli spaces.

An immediate implication of our Main Theorem is that by applying Terasoma's result and using the construction of our Selberg integral-rooted trees functional, we prove that the Hermitian matrix integral as discussed in Mulase [18] evaluates to a Multiple Zeta Value in all 3 cases:

Asymptotically, the Limit as N goes to infinity, and in general.

Furthermore, this construction provides for a positive resolution to Goncharov's conjecture (see [7] pg. 30). The Selberg integral functional can be extended to map the special values L($S^n H^1$(X), n+m) to depth m multiple polylogarithms on X.

**Introduction:**

Our **Main Theorem** is to construct a mapping between the Hopf algebra of Feynman diagrams and the Hopf Algebra of multiple zeta values. This resolves affirmatively an open question whether Feynman diagrams evaluate to multiple zeta values beyond the six loop level. Our work provides positive clarification to this question.

This provides a new functorial environment in which analytical techniques from the Zeta side can be categorically mapped over to the Feynman diagram side and vice-versa. The key step of our theorem is to present the Selberg integral as discussed in Terasoma [22] as a Functional from the Rooted Trees Operad to the Hopf algebra of multiple zeta values. This is a new construction which provides illumination to the relations between zeta values, associators, Feynman diagrams and moduli spaces.

Among other applications of this technique, it provides tools which could be used in the Langlands program, and specifically could be used to construct a "universal" associator over **Q** as discussed in Racinet (see [21] pg. 2). As the applications of Feynman diagrams throughout mathematical physics are many, there are potentially many applications of our new theorem. And conversely, the potential applications of Feynman diagram physics on the zeta function and moduli spaces are without bound.

An immediate implication of our **Main Theorem** is that by applying Terasoma's result and using the construction of our Selberg integral-rooted trees functional, the Hermitian matrix integral as discussed in Mulase [18] evaluates to a multiple zeta value in all 3 cases: Asymptotically, the Limit as N goes to infinity, and in general.

Therefore the 3 methods listed by Mulase are all interrelated and the whole picture requires further investigation. After applying Mellin transforms between the theta and zeta functions and running the logic in reverse, the plan is to further unify the entire three scenarios. This could reveal new understanding of Deligne-Mumford compactifications, Hermitian matrix integrals, and zeta functions in general.

**Background:**

**I. Terasoma's Theorem**

A recent theorem of Terasoma [22] relates Selberg's integral to multiple zeta values by expressing a general form of this integral as an element of the homogeneous ring of MZV's.

Selberg's integral:

$$(I.1) \quad S_\Gamma(V/R, x_i)_{i \in R} = \int_{D(V/R, x_i)_{i \in R}} \Phi(V) \prod_{(i,j) \in E_\Gamma} \alpha_{i,j} \omega_\Gamma$$

Selberg's integral is a function from the free abelian group generated by ordered rooted graphs whose root set and set of vertices are R and V and is denoted by $\Gamma(V,R)$.

Terasoma's proof invokes use of the Drinfeld Associator, which has previously been shown to be expressable in terms of MZV's by Le-Murakami [16].

**II. Mulase's Theorem**

A theorem of Mulase [18], based upon results of Kakei [11] and Goulden-Harer-Jackson [8] expresses the asymptotic expansion of an Hermitian matrix integral in one of three ways:

(II.1) as a Feynman diagram expansion;
(II.2) using classical asymptotic analysis of orthogonal polynomials in the Penner Model resulting in zeta functions;
or (II.3) as a Riemann theta function or a Selberg integral representing a tau-function solution of the KP equations.

**III. Work of Kreimer et.al.**

In [15] Kreimer presents an overview of a broad project involving joint work with Connes, Broadhurst and others focusing on the relationship between the Hopf algebra structure of Feynman diagrams and multi zeta values and knot theory. In [15] he describes the limits to what is known so far regarding evaluating Feynman diagrams in terms of MZV's – that the

Feynman diagrams evaluate up to the 6 Loop level in terms of Euler Zagier Sums, but it is unknown at the 7 Loop level. As of 11/25/01 this was still the known current situation [email from Kreimer 11/25/01].

---

**Our Work**

# MAIN THEOREM

There is an isomorphism between the Hopf Algebra of Feynman diagrams and the Hopf algebra corresponding to the Homogenous Multiple Zeta value ring H in C<<X,Y>> . In other words, Feynman diagrams evaluate to Multi Zeta Values in all cases. This proves a recent conjecture of Connes-Kreimer, and others including Broadhurst and Kontsevich.

**OUTLINE OF PROOF**

Let $\Gamma(V,R)$ be the free abelian group generated by ordered rooted graphs whose root set and set of vertices are R and V.

By **Terasoma's Theorem:** Selberg's integral is a function from $\Gamma(V,R)$ to the Homogenous Multiple Zeta Value ring H.

These rooted graphs $\Gamma(V,R)$ can be given the structure of the Rooted Trees Operad ***RT*** (see Chapoton-Livernet [4]).

Note: The Connes-Kreimer Hopf algebra is dual to the universal enveloping algebra of a Lie algebra generated by the Rooted Trees Operad. See (Chapoton-Livernet [4] pg. 11).

Also note: According to Moch-Uwer-Weinzierl [17] the set of multi-zeta values, which is a subspace of the Euler-Zagier sums, lives in the space of Z-sums. As an aside, the Z-sums also contain the polylogarithms of Goncharov (see M-U-W [17] pg. 6). They show that the coalgebra structure of the Z-sums is identical to the coalgebra structure of the Connes-Kreimer Hopf algebra.

---

**Overview:** First we will use an observation of Mulase (made during the development of his theorem). It states that the Feynman diagram expansion and the Selberg integral expansion result in the same answer for a specific Hermitian Random Matrix model. This implies that the Feynman diagram expansion evaluates to a Selberg integral, and by Terasoma's theorem, subsequently to MZV's. All that remains is to verify that the domain of all the Feynman diagrams in this step is the same as those of our **Main Theorem**. Then apply **Terasoma's theorem**, and we are done.

---------------------------------------------

A **Theorem of Mulase** [18], based upon results of Kakei [11] and Goulden-Harer-Jackson [8] expresses the asymptotic expansion of an Hermitian Matrix integral in one of three ways:

    (M.1) as a Feynman diagram expansion;
    (M.2) using classical asymptotic analysis of orthogonal polynomials in
            the Penner Model resulting in Zeta Functions;
or     (M.3) as a Riemann theta function or a Selberg integral representing a
            tau-function solution of the KP equations.

Note that by an **Observation of Mulase** (pg. 5 in [19]), the evaluation of the matrix integrals agree in the analytical and Feynman diagram case (i.e. cases M.1 and M.2). This provides for necessary and sufficient conditions to infer this covers the full range of diagrams. That is, the domain of step M.1 is **RT**. Similarly, **RT** is also the domain of our **Main Theorem**, as discussed above. Therefore all relevant Feynman diagrams are included in these arguments.

Therefore: The Selberg integral acts as a **functional** from the Operad of Rooted Trees **RT** to the Hopf algebra of MZV's.

Also [aside]: By results of Kontsevich-Soibelman [14] using Strebel differentials; and Kimura-Stasheff-Voronov [13] the Deligne-Knudsen-Mumford moduli space in step M.2 identifies with **RT**. We use this fact to pull-back the Selberg integral functional to the Deligne-Knudsen-Mumford moduli space. [this observation to be used in future work!]

Both $\Gamma(V,R)$ [in Terasoma's paper] and the rooted trees of the Feynman diagrams can be expressed as elements of a suitably defined operad. There exists an automorphism group structure which relates the Feynman diagrams to the MZV's.

Furthermore, investigation into the structure of the Penner model leads to new insight. The ribbon graph evaluation of the Hermitian matrix integral is based upon the topology of the moduli space $M_{g,n} \times R_+^n \cong G_{g,n}^{met}$ which is an orbifold and can be expressed an operad. See Voronov [24]. We extend this mapping from this operad of metric ribbon graphs to the operad of the Feynman diagrams to extend subsequently to MZV's. See Voronov [24] and Kapranov [12].

Furthermore, as this Selberg integral functional maps between the Rooted Trees Operad, and the Connes-Kreimer Hopf algebra, which are equivalent by duality [modulo a twisting operator] we can do more. This functional can then be considered as an **operadic automorphism** on **RT**. The Selberg integral as well as MZV's can then be studied via this [future work]. It remains to be seen what these have to do with the Grothendieck-Teichmuller group. See Bar-Natan [1] and Drinfeld [6].

Note that a Feynman graph with overlapping divergences can be expressed as a sum of Rooted Trees (Connes-Kreimer [5]). In general, the Feynman diagrams form a Hopf Algebra (See Kreimer [15], Connes-Kreimer [5], Varily [23]).

-------------------------------------------

Now we proceed to show Mulase's theorem can be extended to evaluate matrix integrals in terms of multizetas.

-------------------------------------------

**Theorem** [Extension of Mulase's theorem]

The Hermitian matrix integral

1) $\quad Z_n(t,m) =$

$$\int_{H_n} \exp(-\frac{1}{2} Trace(X^2)) \exp(Trace \sum_{j=3}^{2m} \frac{t_j}{j} X^j) \frac{d\mu(X)}{N}$$

is expressable in terms of zeta functions and multiple zeta values.

-------------------------------------------

Sketch of proof:

An outline to the proof of this theorem proceeds as follows:

As Penner's model already evaluates to zeta functions, it remains necessary to do the same for M.1 and M.3. For the Feynman diagram expansion of the Hermitian matrix integral (M.1), experimental results of Broadhurst, Kreimer and Connes-Kreimer evaluate the Feynman diagrams in terms of multi zeta values up to six loops.

Our **MAIN THEOREM** is then used to evaluate Feynman diagrams in terms of multiple zeta values in full generality.

--------------------------------------------

For the tau function solution of the KP equations (M.3) apply Terasoma's theorem to the Selberg integral: It is known that KP equations can be solved in terms of Selberg integrals. Terasoma evaluates Selberg integrals in terms of multiple zeta values.

Mulase's methods (M.2) and (M.3) are essentially equivalent.

That is, the Selberg integral of (M.3) can be related to the zeta functions of (M.2). In addition, results of Broadhurst-Kreimer [3] relate Feynman diagrams in terms of MZV's as well up to 6 loops, and our **MAIN THEOREM** is applied in the general case.

**QED.**

Furthermore, this construction provides for a positive resolution to Goncharov's conjecture (see [7] pg. 30). The Selberg integral functional can be extended to map the special values L($S^n H^1$(X), n+m) to depth m multiple polylogarithms on X.

________________________________________

**Things to do:**

    a) Identify applications in physics with Feynman diagram expansions which provide insight into MZV's and vice versa.

    b) Expand the scope of the **Main Theorem** here to Knot Theory.

    c) Expand the scope of the **Main Theorem** here to pull back the MZV's to provide insight into Deligne Mumford moduli spaces.

    d) Verify the consistency of the **Main Theorem,** and determine the ramifications in our context with the work of Bergere et.al. [2] and Ogreid-Osland [20] and similar papers.

**Acknowledgements**

We extend our appreciation and gratitude to Peter Bouwnegt (Univ of Adelaide), Pierre Deligne (IAS), Vladimir Drinfeld (Univ of Chicago), Motohico Mulase (Univ of Calif./Davis) and Donald Richards (IAS/Univ of Va.) for discussions and lectures relevant to this paper.

In particular, Donald conjectured to me in the Spring of 2001 about a stronger version of Terasoma's theorem in which Selberg integrals are evaluated explicitly in terms of MZV's.  It was this conjecture which led to my deeper investigations into Terasoma's paper. The conjecture is still part of my work in progress. Deligne gave several lectures on Racinet's paper and we had several discussions on MZV's.  Mulase explained to me the unpublished details on the proof of his theorem based upon work of Kakei and Goulden-Harer-Jackson,  Bouwnegt gave me some advice on the KZ equations and field theory, and Drinfeld gave me a very helpful lecture on associators.

I would also like to express thanks to Dennis Sullivan for hosting his series of guest lectures on Quantum Algebra, wherein I was introduced to many of the topics mentioned in this note.

And finally tremendous thanks to Michael Anshel for everything, including his advice and general direction to look towards zeta functions as a unifying force.

## References

1) D.Bar-Natan, On Associators and the Grothendieck-Teichmuller Group I, Q-Alg/9606021, 7/1/96.

2) M.C.Bergere, C.deCalan and A.P.C.Malbouisson, A Theorem on Asymptotic Expansion of Feynman Amplitudes, Commun. Math. Phys 62, 137-158, (1978).

3) D.J.Broadhurst and D.Kreimer,  Association of Multiple Zeta Values with Positive Knots via Feynman Diagrams up to 9 Loops,  Hep-th/9609128, 11/18/96.

4) F.Chapoton and M.Livernet, Pre-Lie Algebras and the Rooted Trees Operad, QA/0002069, 6/29/00.

5) A.Connes and D.Kreimer,  Hopf Algebras, Renormalization and Noncommutative Geometry, 9808042, 7/7/98.

6) V.Drinfeld, On Quasitriangular Quasi-Hopf algebras and a group closely connected with Gal(Q/Q), Leningrad Math. J. Vol 2 (1991), No. 4, 829-860.

7) A.B.Goncharov, Multiple Zeta Values, Galois Groups, and Geometry of
  Modular Varieties, AG/0005069, 9/14/00.

8) I.P.Goulden, J.L.Harer, D.M.Jackson, A Geometric Parametrization for
  Virtual Euler Characteristic … , AG/9902044, 2/5/99.

9) M.E.Hoffman, Algebras of Multiple Zeta Values, Quasi-Symmetric Functions,
  and Euler Sums, Talk at UQAM US Naval Academy, 5/1/98.
  [available from his webpage].

10) M.E.Hoffman and Y.Ohno, Relations of Multiple Zeta Values and their
  Algebraic Expression, QA/0010140, 10/13/00.

11) S.Kakei Orthogonal and Symplectic Matrix Integrals and Coupled KP
  Heirarchy, Solv-Int/9909023, 9/22/99.

12) M.Kapranov, Operads and Algebraic Geometry, Doc. Math. J. Extra Vol.
  ICM 1998 II, 277-286.

13) T.Kimura, J.Stasheff and A.A.Voronov, On Operad Structures of Moduli
  Spaces and String Theory, Hep-Th/9307114, 3/14/94.

14) M.Kontsevich and Y.Soibelman, Deformations of Algebras over Operads
  and Deligne's Conjecture, QA/0001151, 1/27/00.

15) D.Kreimer, Combinatorics of (perturbative) Quantum Field Theory,
  Hep-th/0010059, 10/9/00.

16) T.T.Q.Le and J.Murakami, Kontsevich's integral for the Kauffman polynomial,
  Nagoya J. Math 142 (1996), 39-65.   See also H.Furusho, The Multiple Zeta
  Value Algebra …, NT/0011261 pg. 9 about an error in Le-Murakami.

17) [M-U-W] S.Moch, P.Uwer and S.Weinzierl:  Nested Sums, Expansion of
  Transcendental Functions and Multi-Scale Multi-Loop Integrals,
  Hep-Th/0110083, 10/5/01.

18) M.Mulase, Lectures on the Asymptotic Expansion of a Hermitian
  Matrix Integral, MathPhys/9811023, 11/25/98.

19) M.Mulase, Lectures on Combinatorial Structure of Moduli Spaces of
  Riemann Surfaces,  UC Davis CA Preprint, [Work in Progress] 11/1/01.

20) O.M.Ogreid and P.Osland, Some Infinite Series related to Feynman
  Diagrams, Math-Ph/0010026, 10/20/00.

21) G.Racinet, Torseurs associes a certaines relations algebriques entre

polyzetas aux racines de l'unite, QA/0012024, 3/5/01.

22) T.Terasoma, Selberg Integral and Multiple Zeta Values, AG/9908045, 3/25/01.

23) J.C.Varilly, Hopf Algebras in Noncommutative Geometry, Hep-Th/0109077, 9/10/01.

24) A.A.Voronov, Notes on Universal Algebra, QA/0111009, 11/1/01.

David H. Wohl

Hadar HaTorah Yeshiva
824 Eastern Parkway
Brooklyn, NY
11213

Department of Mathematics
Touro College
1602 Avenue J
Brooklyn, NY 11230

DAVIDHW@TOURO.EDU

DHWSYSTEMS@ATT.NET